\documentclass[12pt]{amsart}
\font\emailfont=cmtt10

\usepackage{amsmath,amsthm,amsfonts,amscd,flafter,epsf}

\newcommand\commentable[1]{#1}

\newcommand{\rk}{\mathrm{rk}}

\newtheorem{theorem}{Theorem}[section]
\newtheorem{prop}[theorem]{Proposition}
\newtheorem{cor}[theorem]{Corollary}
\newtheorem{conj}[theorem]{Conjecture}
\newtheorem{lemma}[theorem]{Lemma}

\newtheorem{defn}[theorem]{Definition}

\def\endproof{\relax\ifmmode\expandafter\endproofmath\else
  \unskip\nobreak\hfil\penalty50\hskip.75em\hbox{}\nobreak\hfil\bull
  {\parfillskip=0pt \finalhyphendemerits=0 \bigbreak}\fi}
\def\endproofmath$${\eqno\bull$$\bigbreak}
\def\bull{\vbox{\hrule\hbox{\vrule\kern3pt\vbox{\kern6pt}\kern3pt\vrule}\hrule}}

\newcounter{bean}
\newcommand{\Q}{\mathbb{Q}}
\newcommand{\R}{\mathbb{R}}

\newcommand{\Z}{\mathbb{Z}}

\newcommand{\Zmod}[1]{\Z/{#1}\Z}

\newcommand{\cm}{\cdot}

\newcommand{\Nbd}[1]{{\mathrm{nd}}(#1)}
\newcommand{\nbd}[1]{\Nbd{#1}}

\newcommand{\ModSWfour}{\mathcal{M}}
\newcommand{\ModFlow}{\ModSWfour}

\newcommand{\SpinC}{{\mathrm{Spin}}^c}

\newcommand\abuts\Rightarrow
\newcommand\Sym{\mathrm{Sym}}

\newcommand\HFpRed{\HFp_{\red}}

\newcommand\relspinc{\underline{\spinc}}

\newcommand\Filt{\mathcal F}

\newcommand\ModSphere{\ModFlow\left({\mathbb S}\longrightarrow 
\Sym^{g-1}(\Sigma_{1})\times \Sym^2(\Sigma_{2})\right)}
\newcommand\ModSpheres\ModSphere

\newcommand\CFa{\widehat{CF}}

\newcommand{\red}{\mathrm{red}}

\newcommand\HFp{\HFb}

\newcommand\HFa{\widehat{HF}}
\newcommand\HFb{HF^+}

\newcommand\UnparModSp{\widehat \ModSp}
\newcommand\UnparModFlow\UnparModSp
\newcommand\Mod\ModSp

\newcommand\PD{\mathrm{PD}}

\newcommand{\spinc}{\mathfrak s}

\newcommand\ModMaps{\mathcal M}
\newcommand\ModSp\ModMaps

\newcommand\Field{\mathbb F}

\newcommand\Dual{\mathcal D}
\newcommand\Duality\Dual

\newcommand\orb{\mathrm{orb}}

\newcommand\spincrel\relspinc

\newcommand\CFK{CFK}
\newcommand\HFK{HFK}

\newcommand\CFKa{\widehat\CFK}

\newcommand\CFKinf{\CFK^{\infty}}

\newcommand\HFKa{\widehat\HFK}

\commentable{

\title[{On knot Floer homology and lens space surgeries}]
{On knot Floer homology and lens space surgeries}

\author[Peter Ozsv{\'a}th]{Peter Ozsv\'ath}
\address{Department of
Mathematics, Columbia University, New York 10027 \newline
\indent{\emailfont{petero@math.columbia.edu}}}

\author[Zolt{\'a}n Szab{\'o}]{Zolt{\'a}n Szab{\'o}} 
\address{Department of
Mathematics, Princeton University, New Jersey 08540 \newline
\indent{\emailfont{szabo@math.princeton.edu}}}}

\begin{document}

\begin{abstract}  
In an earlier paper, we used the absolute grading on Heegaard Floer
homology $\HFp$ to give restrictions on knots in $S^3$ which
admit lens space surgeries. The aim of the present article is to
exhibit stronger restrictions on such knots, arising from knot
Floer homology. One consequence is that the non-zero coefficients
of the Alexander polynomial of such a knot are $\pm 1$.  This
information can
in turn
be used to prove that certain lens spaces are
not obtained as integral surgeries on knots. In fact, combining our
results with constructions of Berge, we classify lens spaces $L(p,q)$
which arise as integral surgeries on knots in $S^3$ with $|p|\leq
1500$. Other applications include bounds on the four-ball genera of
knots admitting lens space surgeries (which are sharp for Berge's
knots), and a constraint on three-manifolds obtained as integer
surgeries on alternating knots, which is closely to related to a
theorem of Delman and Roberts.
\end{abstract}

\maketitle
\section{Introduction}

Let $K\subset S^3$ be a knot for which some integral surgery gives a
lens space $L(p,q)$. The surgery long exact sequence for 
Heegaard Floer homology $\HFp$,
together with the absolute grading on the latter group, can be
combined to give a number of restrictions on $K$,
see~\cite{AbsGraded}.

Consequently, for each fixed lens space $L(p,q)$, there is an
explicitly determined, finite list of symmetric polynomials which
might arise as the Alexander polynomials of such knots. The
indeterminacy can be clarified from the following dual point of
view. A knot in $S^3$ whose surgery gives $L(p,q)$ induces a knot $K'$
in $L(p,q)$ on which some surgery gives $S^3$. The stated
indeterminacy, then, corresponds to the possible homology classes for
$[K']\in H_1(L(p,q);\Z)$.  Indeed, there are straightforward
homological obstructions to realizing a given homology class in
$H_1(L(p,q);\Z)$ in this way from
a knot in $S^3$ (or indeed from any integer homology
three-sphere).

The results of~\cite{AbsGraded} go beyond these homological
considerations to give additional constraints on the Alexander
polynomials of the knots $K$. These further constraints are specific
to $S^3$: they can be used to rule out lens space surgeries even in
cases where the lens space is realized as a surgery on a knot in some
other integral homology three-sphere.

The aim of the present article is to strengthen considerably these
constraints, with the help of the Floer homology of knots, see~\cite{Knots}
and also~\cite{RasmussenThesis}. In terms of
the Alexander polynomial, our results here show that if $K$ is a knot
with the above properties, then all the coefficients of its Alexander
polynomial are $\pm 1$, and the non-zero coefficients alternate in
sign. Actually, since our results apply to a
wider class of three-manifolds than lens spaces, before stating the
theorems precisely, we discuss the class of three-manifolds we study.

\subsection{Knot Floer homology and $L$-space surgeries.} 
The appropriate generalization of the notion of lens spaces, for our
purposes, is given in the following definition.  Note that $\HFa(Y)$
is the three-manifold invariant defined in~\cite{HolDisk}.

\begin{defn}
\label{def:LSpaces}
A closed three-manifold $Y$ is called an $L$-space if $H_1(Y;\Q)=0$,
and $\HFa(Y)$ is a free Abelian group
whose rank coincides with the number of elements in
$H_1(Y;\Z)$, which we write as $|H_1(Y;\Z)|$.
\end{defn}

The set of $L$-spaces includes all lens spaces $L(p,q)$ and, indeed,
all spaces with ``elliptic geometry,'' i.e. all the finite, free
quotients of $S^3$ by groups of isometries
(c.f. Proposition~\ref{prop:SpaceForms} below). It also includes a
class of plumbing manifolds which are obtained as plumbings specified
by trees, for which the surgery coefficient associated to each vertex
is no smaller than the number of edges meeting at that vertex
(according to Theorem~\ref{HFSymp:thm:FloerHomology}
of~\cite{HolDiskSymp}). The set of $L$-spaces is closed under
connected sums, and the following additional operation: fix an
$L$-space $Y$, and a knot $K\subset Y$ with a choice of framing
$\lambda$ for which
$$|H_1(Y_{\lambda+\mu}(K))|=|H_1(Y)|+|H_1(Y_{\lambda}(K))|,$$ where
$\mu$ denotes the meridian for the knot, and $Y_{\lambda}(K)$ denotes
the three-manifold obtained from $Y$ by performing surgery on $Y$
along $K$ with framing $\lambda$.  Then, if both $Y$ and
$Y_{\lambda}(K)$ are $L$-spaces, then so is $Y_{\lambda+\mu}(K)$.
This construction gives infinitely many hyperbolic $L$-spaces. For
instance, let $P(a,b,c)$ denote the three-stranded pretzel knot with
$a$, $b$, and $c$ twists respectively.  As observed by Fintushel and
Stern~\cite{FSLensSurgeries}, $+18$ surgery on $P(-2,3,7)$ is a lens
space.  Thus, applying the above principle to the knot $K$ in the
$L$-space $S^3$, and induction, we see that for all integers $p\geq
18$, $S^3_p(P(-2,3,7))$ is an $L$-space. These are hyperbolic for all
sufficiently large $p$, according to a theorem of Thurston~\cite{Thurston},
\cite{Thurston2}
(in fact, the fundamental group is infinite for all $p>19$,
c.f.~\cite{Mattman}).  A more in-depth discussion of $L$-spaces with
more examples is given in Section~\ref{sec:LSpaces}.

The results of this paper are built on the following theorem about
the Floer homology of a knot which admits an $L$-space surgery.
To state the result, recall that 
there is a knot Floer homology group associated to 
a knot $K$ in $S^3$ and an integer $i$, 
which is a graded Abelian group, denoted
$\HFKa(K,i)$, c.f.~\cite{Knots}, see
also~\cite{RasmussenThesis}.

\begin{theorem}
\label{thm:FloerHomology}
Suppose that $K\subset S^3$ is a knot for which there is a positive
integer $p$ for which $S^3_p(K)$ is an $L$-space. Then, there is an
increasing sequence of non-negative integers
$$n_{-k}<...<n_k$$
with the property that $n_i=-n_{-i}$, with the following significance.
If for $-k\leq i \leq k$ we let
$$\delta_{i}=\left\{\begin{array}{ll}
0 & {\text{if $i=k$}} \\
\delta_{i+1}-2(n_{i+1}-n_{i})+1 &{\text{if $k-i$ is odd}} \\
\delta_{i+1}-1 & {\text{if $k-i>0$ is even,}}
\end{array}\right.$$
then
$\HFKa(K,j)=0$ unless $j=n_i$ for some $i$, in which case
$\HFKa(K,j)\cong \Z$ and it is supported entirely in dimension $\delta_i$.
\end{theorem}

Since 
$$\sum_i\chi(\HFKa(K,i)) \cm T^i=\Delta_K(T)$$ is the symmetrized Alexander polynomial
(c.f. Proposition~\ref{Knots:prop:Euler} of~\cite{Knots}),
the above theorem says that $\HFKa$ is determined explicitly from the Alexander polynomial
of $K$. Conversely, the above theorem gives strong restrictions on the
Alexander polynomials of knots which
admit $L$-space surgeries:

\begin{cor}
\label{cor:StructAlex}
Let $K\subset S^3$ be a knot for which there is an integer $p$ for which
$S^3_p(K)$ is an $L$-space. Then the Alexander polynomial of $K$ has the form
$$\Delta_K(T) = (-1)^k+ \sum_{j=1}^k(-1)^{k-j} (T^{n_j}+T^{-n_j}),$$
for some increasing sequence of positive integers $0<n_1<n_2<...<n_k$.
\end{cor}

For a fixed $L$-space $Y$, the possible polynomials which could occur
as the Alexander polynomials of knots $K\subset S^3$ for which
$S^3_p(K)\cong Y$ is determined up to a finite indeterminacy by the
absolute grading on $\HFa(Y)$ (c.f.~\cite{AbsGraded}, but observe that
this result also follows from the methods of the present paper, see
Section~\ref{sec:Proof}).
Thus, the above corollary can be used to give new restrictions on
which $L$-spaces arise as $+p$ surgeries on knots in $S^3$.

\subsection{Alexander polynomials and lens space surgeries}

As an illustration, let $d(L(p,q),i)$ denote the absolute grading of
the generator of $\HFa(L(p,q),i)$.  (Here, we use the orientation
convention that $L(p,q)$ is obtained by $p/q$ surgery on the unknot in
$S^3$.)  We showed in~\cite{AbsGraded}
(c.f. Proposition~\ref{AbsGraded:prop:dLens}; compare
also~\cite{MilnorTorsion} and~\cite{Turaev}), that this quantity is
determined by the recursive formula
\begin{eqnarray*}
d(-L(1,1),0)&=&0\\
d(-L(p,q),i)&=&\left(\frac{pq-(2i+1-p-q)^2}{4pq}\right)-d(-L(q,r),j), 
\end{eqnarray*}
where $r$ and $j$ are the reductions modulo $q$ of $p$ and $i$
respectively. Note that we are implicitly using here a specific
identification $\Zmod{p}\cong \SpinC(L(p,q))$ (defined explicitly in
Subsection~\ref{AbsGraded:subsec:dLens} of~\cite{AbsGraded}, but not
crucial for our purposes here). We have the following consequence of
Corollary~\ref{cor:StructAlex}:

\begin{cor}
\label{cor:LensCondition}
The lens space $L(p,q)$ is obtained as surgery on a knot $K\subset S^3$
only if there is a one-to-one correspondence 
$$\sigma\colon \Zmod{p}\longrightarrow \SpinC(L(p,q))$$
with the following symmetries:
\begin{itemize}
\item $\sigma(-[i])={\overline{\sigma([i])}}$
\item there is an isomorphism $\phi\colon \Zmod{p}
\longrightarrow \Zmod{p}$
with the property that 
$$\sigma([i])-\sigma([j])=\phi([i-j]),$$
\end{itemize}
with the following properties.  For $i\in\Z$, let $[i]$ denote its reduction 
modulo $p$,
and define
$$t_i=
\left\{\begin{array}{ll}
-d(L(p,q),\sigma[i]) + d(L(p,1),[i]) &{\text{if $2|i|\leq p$}} \\ \\
0 & {\text{otherwise,}}
\end{array}
\right. 
$$ then the Laurent polynomial
$$1 + \sum_{i} \left(\frac{t_{i-1}}{2}-t_i+ \frac{t_{i+1}}{2}\right)
T^i =\sum_i a_i \cm T^i$$ has integral coefficients,
all of which satisfy $|a_i|\leq 1$, and all of
its non-zero coefficients alternate in sign.
\end{cor}

For instance, a straightforward if tedious calculation using the above
result shows that $L(17,2)$ does not occur as integral surgery on any
knot in $S^3$, even though it passes all the criteria
from~\cite{AbsGraded} (in particular, it is realizable as $+17$
surgery on a knot in some other integral homology three-sphere).
Similar remarks hold for $L(19,17)$ and $L(26,23)$ (compare the list
at the end of~\cite{AbsGraded}).

In fact, these obstructions are particularly powerful when one
combines them with Berge's construction of knots which admit lens
space surgeries, see~\cite{Berge}. Indeed, there is evidence
suggesting that the conditions on $L(p,q)$ in
Corollary~\ref{cor:LensCondition} necessary for it to be realized as
integral surgery on a knot in $S^3$ are also sufficient.  We return to
this point at the end of the present introduction, after describing
Berge's constructions. But first, we turn to some other immediate
applications of Theorem~\ref{thm:FloerHomology}.

\subsection{Alternating knots and $L$-space surgeries.}
In another direction, we obtain the following consequence of
Corollary~\ref{cor:StructAlex} (together with properties of the
Alexander polynomials of alternating knots,
c.f. Section~\ref{sec:AltKnots}), which is rather similar in spirit to
a theorem of Delman and Roberts~\cite{DelmanRoberts} obtained using
the theory of laminations:

\begin{theorem}
\label{thm:AltLens}
If $K\subset S^3$ is an alternating knot with the property
that
some integral 
surgery along $K$ is an
$L$-space, then $K$ is a $(2,2n+1)$ torus knot, for some
integer $n$.
\end{theorem}

\subsection{Bounding the four-ball genus}
To go beyond the Alexander polynomial, recall that the knot Floer
homology $\bigoplus_m \HFKa(K,m)$ is the homology of the graded
complex associated to a filtration $$...\subseteq \Filt(K,m)\subseteq
\Filt(K,m+1) \subseteq ...$$ of the chain complex $\CFa(S^3)$ (whose
homology is $\Z$, in a single dimension), induced by the knot
$K$. This filtration gives an integer $\tau(K)$ which is defined to be
the smallest integer $m$ for which the induced map on homology
$$\iota^m_K\colon H_*(\Filt(K,m))\longrightarrow \HFa(S^3)\cong \Z$$ is
non-trivial. In~\cite{4BallGenus} 
(c.f. Corollary~\ref{4BallGenus:cor:GenusBounds}
of~\cite{4BallGenus}) and also~\cite{RasmussenThesis}, 
it is shown  that if $g^*(K)$ denotes the
four-ball genus, then
\begin{equation}
\label{eq:FourBallGenusBound}
|\tau(K)|\leq g^*(K).
\end{equation}
Note that $g^*(K)$ gives a lower bound on the unknotting number of $K$.
Combining Inequality~\eqref{eq:FourBallGenusBound}
with Theorem~\ref{thm:FloerHomology}, we obtain the following:

\begin{cor}
\label{cor:CalcTau}
Suppose that $K\subset S^3$ is a knot which admits an integral
$L$-space surgery, then $|\tau(K)|$ coincides with the degree of the
symmetrized Alexander polynomial of $K$.  In particular, the four-ball
genus $g^*(K)$ is bounded below by this degree.
\end{cor}

\begin{proof}
The first claim follows immediately from the description of the knot Floer
homology given in Theorem~\ref{thm:FloerHomology}, while the second follows
from the first, together with Inequality~\eqref{eq:FourBallGenusBound}.
\end{proof}

\begin{figure}
\mbox{\vbox{\epsfbox{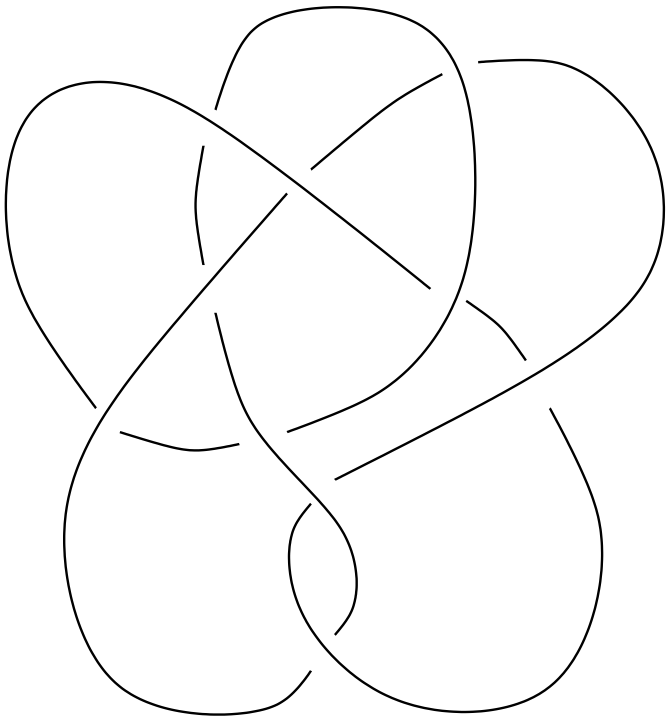}}}
\caption{\label{fig:Ten132}
{\bf{The knot $10_{132}$.}} This knot has Alexander polynomial
$1-(T+T^{-1})+(T^2+T^{-2})$, and unknotting number one.}
\end{figure}

Corollary~\ref{cor:CalcTau} also gives an illustration of how
Theorem~\ref{thm:FloerHomology} goes beyond the Alexander
polynomial. As an amusing application, consider the knot $K=10_{132}$
pictured in Figure~\ref{fig:Ten132}, the ten-crossing knot whose
Alexander polynomial is $$\Delta_K(T)=T^{-2}-T^{-1}+1-T+T^2.$$ This
Alexander polynomial satisfies the criteria of
Corollary~\ref{cor:StructAlex}. However, the knot clearly has
unknotting number one, and hence according to
Corollary~\ref{cor:CalcTau}, this knot admits no $L$-space surgeries.

Let $T_{p,q}$ denote the $(p,q)$ torus knot.  Since $pq\pm 1$-surgery
on the $(p,q)$ torus knot is a lens space, the above corollary shows
that $\tau(T_{p,q})=\frac{(p-1)(q-1)}{2}$, and hence (after a careful
choice of unknotting) that the unknotting number of $T_{p,q}$ is given
by this quantity. This result was first proved by Kronheimer and
Mrowka~\cite{KMMilnor} (and conjectured by Milnor~\cite{Milnor}).

In fact, Corollary~\ref{cor:CalcTau} gives a calculation of the
four-ball genera of all knots coming from Berge's constructions,
see~\cite{Berge}.  Specifically, recall that Berge's constructions
have a particularly nice description from the point of view of knots
in lens spaces.

\begin{defn} 
\label{def:BergeKnots}
Consider the standard genus one Heegaard diagram for $L(p,q)$, where
the two attaching circles $\alpha$ and $\beta$ meet in exactly $p$
points. A {\em lens space Berge knot} $K'\subset L(p,q)$ is one which
is formed from a pair of arcs, one of which is supported in the
attaching disk for $\alpha$ and the other is supported in the
attaching disk for $\beta$, with the additional property that $[K']\in
H_1(L(p,q);\Z)\cong \Zmod{p}$ is a generator.
\end{defn}

The following result is verified in Section~\ref{sec:BergeKnots},
using results of Stallings~\cite{Stallings} and Brown~\cite{Brown}:

\begin{prop}
\label{prop:BergesFibered}
All lens space Berge knots are fibered.
\end{prop}

\begin{defn}
When integral surgery of $L(p,q)$ along some lens space Berge knot
$K'$ gives $S^3$, there is a naturally induced knot $K\subset S^3$ for
which some integral surgery gives $L(p,q)$. We call this induced knot
a {\em classical Berge knot}.
\end{defn}

\begin{cor}
\label{cor:Berge4BallGenus}
Let $K$ be a classical Berge knot. Then, the degree of the Alexander
polynomial agrees with both the Seifert and four-ball genera of $K$.
\end{cor}

\begin{proof} Since $K$ is fibered, its Seifert genus agrees with
the degree of its Alexander polynomial. The statement about the four-ball
genus now follows from Corollary~\ref{cor:CalcTau}.
\end{proof}

\subsection{Realizing lens spaces}

In~\cite{Berge}, Berge proves the
following theorem:

\begin{theorem}(Berge)
The lens space $L(p,q)$ arises as integral surgery on a knot in $S^3$ if 
we can find integers
$A,a,B,b$ so that $p=|Aa+Bb|$, $a^2 q\equiv \pm b^{\pm 2}\pmod{p}$,
which satisfy at least one of the following additional constraints
\begin{list}
{(\arabic{bean})}{\usecounter{bean}\setlength{\rightmargin}{\leftmargin}}
\item $A=1$, $a=\pm 1$, $(B,b)=1$, $B\geq 2$,
\item $A=1$, $a=\pm 1$, $(B,b)=2$, $B\geq 4$,
\item  $A>1$, $a=\pm 1$,  and there is an integer $\epsilon=\pm 1$ so that
$(B+\epsilon)/A$ is an odd integer and $b\equiv -2\epsilon Aa\pmod{B}$,
\item $A>3$, $a=\pm 1$ and there is an $\epsilon=\pm 1$ so that $(2B+\epsilon)/A$ is 
integral, and $b\equiv -\epsilon A a\pmod{B}$,
\item $A>1$, $A$ is odd, $a=\pm 1$, and there is an $\epsilon=\pm 1$ such that
$(B-\epsilon)/A$ is an integer, and $b\equiv -\epsilon A a \pmod{B}$,
\item $A>2$, $A$ is even, $a=\pm 1$ $B=2A+1$, and $b\equiv -a(A-1)\pmod{B}$,
\item $a=-(A+B)$, $b=-B$,
\item $a=-(A+B)$, $b=B$,
\item $(A,B,a,b)=(4J+1,2J+1,6J+1,-J)$ for some integer $J$,
\item $(A,B,a,b)=(6J+2,2J+1,4J+1,-J)$ for some integer $J$,
\item $(A,B,a,b)=(6J+4,2J+1,-4J-3,J+1)$ for some integer $J$,
\item $(A,B,a,b)=(4J+3,2J+1,-6J-5,J+1)$ for some integer $J$.
\end{list}
\end{theorem}

Berge proves the above theorem by explicitly constructing the
corresponding knots in $S^3$. For instance, lens spaces of Type~(1)
are realized by surgeries on torus knots, Type~(2) by surgeries on cables of
torus knots, Types~(3)-(6) by other knots supported in a solid torus
(for which some surgery gives another solid torus), Type~(7) by surgeries on 
knots supported in the Seifert surface of a trefoil, Type~(8) by surgeries
on knots supported in the Seifert surface of the figure eight knot, and
Types~(9)-(12) are some other ``sporadic'' examples.
Experimental verification suggests the following purely
combinatorial conjecture (which we have verified for $|p|\leq 1500$
using a program written in Mathematica~\cite{Mathematica}):

\begin{conj}
A lens space $L(p,q)$ appears on Berge's list above if and only
if it passes the conditions of Corollary~\ref{cor:LensCondition}.
\end{conj}

A proof of the above conjecture, of course, would prove that
Berge's conditions on a lens space are necessary and sufficient for it
to be realized as integral surgery on a knot in $S^3$. Thus, our computer
verification can be alternately phrased as follows:

\begin{prop}
\label{prop:Experiment}
For all $p\leq 1500$, the lens spaces which are realized as integer
surgeries on knots in $S^3$ are precisely those lens spaces which are
realized on Berge's list.
\end{prop}

Note that Berge conjectures a stronger statement: he conjectures that
his scheme~\cite{Berge} enumerates all knots which admit lens space
surgeries. The above proposition could be viewed as partial evidence
supporting his conjecture.

\subsection{Further remarks}

Recall that in~\cite{Contact}, we
proved that if $K$ is a fibered knot with genus $g$, then
$\HFKa(K,g)\cong\Z$, generalizing the standard fact that a fibered
knot has monic Alexander polynomial.  Thus
Theorem~\ref{thm:FloerHomology} could be seen as evidence supporting
the conjecture that all knots with lens space surgeries are fibered.
(Note that at the time of the writing of this paper, the authors know
of no non-fibered knots for which $\HFKa(K,d)\cong \Z$, where $d$
denotes the degree of the Alexander polynomial of $K$.)

\subsection{Organization.} We discuss  $L$-spaces in Section~\ref{sec:LSpaces}.
In Section~\ref{sec:Proof} we prove Theorem~\ref{thm:FloerHomology}
and its immediate corollaries, Corollaries~\ref{cor:StructAlex}
and~\ref{cor:LensCondition}, and also Proposition~\ref{prop:Experiment}.
In Section~\ref{sec:AltKnots}, we prove Theorem~\ref{thm:AltLens}
from Theorem~\ref{thm:FloerHomology}, and a result on the Alexander
polynomials of alternating knots. Finally, in Section~\ref{sec:BergeKnots},
we verify Proposition~\ref{prop:BergesFibered}.

\subsection{Acknowledgements.} We would like to thank Kenneth Baker,
Andrew Casson, Danny Calegary, David Gabai, Cameron Gordon, John
Luecke, Paul Melvin, and Walter Neumann for interesting discussions
during the course of this work.  In particular, it was Gabai who
called to our attention the fact that all known knots giving lens
space surgeries are fibered; we are also indebted to Calegari and
Neumann for explaining to us the work of Kenneth Brown~\cite{Brown}
used in Section~\ref{sec:BergeKnots}.

\section{$L$-spaces}
\label{sec:LSpaces}

The aim of the present section is to collect some of the key
properties of $L$-spaces, and to
give some constructions. In fact, much of the material here is not
new, but can be found sprinkled throughout most of our papers on
Heegaard Floer homology. It is for this reason that we feel that it
might be useful to collect the properties in one place.

Recall that $\HFa(Y)$ is a finitely generated, $\Zmod{2}$-graded
Abelian group which, for rational
homology three-spheres, satisfies the relation that
\begin{equation}
\label{eq:EulerHFa}
\chi(\HFa(Y))=|H_1(Y;\Z)|
\end{equation}
(c.f. Proposition~\ref{HolDiskTwo:prop:EulerHFa} of~\cite{HolDiskTwo})
and, in particular, for any rational homology three-sphere,
$|H_1(Y;\Z)|\leq \rk \HFa(Y)$.  When $Y$ is an $L$-space, this
inequality is an equality.

In general, $\HFa(Y)$ depends on the orientation used for $Y$, but its
total rank does not (c.f. Proposition~\ref{HolDiskTwo:prop:Duality}
of~\cite{HolDiskTwo}). Using coefficients in an arbitrary field, we
see that the condition of being an $L$-space is also independent of
the orientation of $Y$. Note that $L$-spaces have an alternate
characterization in terms of other elements of the Heegaard Floer
homology package defined in~\cite{HolDisk}: a rational homology
three-sphere is an $L$-space if and only if $\HFpRed(Y)=0$. We will
have no further use for this characterization in the present paper,
but point it out as it appears to be the characterization which
generalizes more neatly to the case where $b_1(Y)>0$.

The fact that lens spaces are $L$-spaces follows immediately
from
their standard genus one Heegaard diagrams
(c.f. Proposition~\ref{HolDiskTwo:prop:Lensspaces}
of~\cite{HolDiskTwo}).

The functor $\HFa(Y)$ enjoys a K{\"u}nneth principle for connected
sums (Proposition~\ref{HolDiskTwo:prop:ConnSum} of~\cite{HolDiskTwo}),
from which it follows readily that the set of
$L$-spaces is closed under connected sums.

Suppose that $K\subset Y$ is a knot in a rational homology
three-sphere, and let $\mu$ be the meridian for $K$ and let $\lambda$
be any choice of longitude (i.e. simple, closed curve in the torus
$\partial \nbd{K}$ which meets $\mu$ in a single transverse point of
intersection). Indeed, suppose that $\lambda$ is chosen so that the
three-manifolds $Y_\lambda(K)$ and $Y_{\lambda+\mu}(K)$ are both
rational homology three-spheres.  We have the following result
(compare Lemma~\ref{AbsGraded:lemma:InvisQSphere}
of~\cite{AbsGraded}):

\begin{prop}
\label{prop:LSpaces}
Let $K\subset Y$ be a knot in a rational homology three-sphere,
and let $\lambda$ be a choice of longitude for the knot, so that 
$Y_{\lambda}(K)$, $Y_{\lambda+\mu}(K)$ are also rational homology three-spheres,
and 
$$|H_1(Y_{\lambda+\mu}(K))|=|H_1(Y)|+|H_1(Y_{\lambda}(K))|.$$
If $Y$ and $Y_{\lambda}(K)$ are $L$-spaces, then so is
$Y_{\lambda+\mu}(K)$.
\end{prop}

\begin{proof}
This follows readily from the long exact surgery sequence for $\HFa$
(c.f. Theorem~\ref{HolDiskTwo:thm:GeneralSurgeryHFa}
of~\cite{HolDiskTwo}), which in the present case reads: $$
\begin{CD}
...@>>>\HFa(Y) @>>> \HFa(Y_{\lambda}(K))@>>>\HFa(Y_{\lambda+\mu}(K))@>>>...
\end{CD}
$$

In particular, we see that $$\rk\HFa(Y_{\lambda+\mu}(K))\leq
\rk\HFa(Y)+\rk\HFa(Y_{\lambda}(K)).$$ It follows immediately that if
$Y$ and $Y_{\lambda}(K)$ are $L$-spaces, then
$$\rk\HFa(Y_{\lambda+\mu}(K))\leq |H_1(Y_{\lambda+\mu}(K))|,$$ while the
opposite inequality is provided by Equation~\eqref{eq:EulerHFa},
forcing equality to hold. Moreover,  repeating the above
argument with coefficients in any finite field,
one verifies that $\HFa(Y_{\lambda+\mu}(K))$ has no torsion.
\end{proof}

The above proposition guarantees that if $K\subset S^3$ is a knot
with the property that 
$S^3_p(K)$ is an $L$-space,
for some positive
integer $p$, then so is $S^3_n(K)$ for all integers $n\geq
p$.

In~\cite{SomePlumbs}, we give a characterization of Seifert
fibered $L$-space which we recall presently.
Recall that a
Seifert fibered rational homology three-sphere is specified by a
collection of integers $b$, $\{\alpha_i\}_{i=1}^n$, and $\{\beta_i\}_{i=1}^n$
(sometimes abbreviated $(b; \beta_1/\alpha_1, ...,\beta_n/\alpha_n)$),
where here all $\alpha_i>2$ and the $0<\beta_i<\alpha_i$, and
$(\alpha_i,\beta_i)=1$ (see~\cite{Seifert}, see also~\cite{Scott} for
a modern treatment). The $\{\alpha_i\}$ specify the base orbifold
(which in the present case must have genus zero). The number $n$
is the number of singular orbits for the circle action on $Y$.

A Seifert fibered space has an orbifold degree given by the formula
$$b+\sum_{i}\frac{\beta_i}{\alpha_i}.$$ When the base has genus zero,
the orbifold degree is non-zero precisely when $Y$ is a rational
homology three-sphere. Note that the orbifold degree changes sign
under orientation reversal of $Y$.

A Seifert space $(b; \beta_1/\alpha_1,...\beta_n/\alpha_n)$ can be
realized as the boundary of a plumbing of spheres, where the spheres
are arranged in a star-like pattern, so that the central node has
self-intersection number $b$, and the multiplicities of the chains of
spheres is given by the Hirzebruch-Jung fractional expansion of
$\alpha_i/\beta_i$. Let $G$ denote this labeled graph: i.e. this is a
tree equipped with a function $m$ from the vertices of $G$ to the
integer, which gives rise in the usual manner to an inner product on
the vector space $V$ generated by the vertices.  In topological terms,
$V$ is $H_2(W(G))$ where $W(G)$ is the four-manifold constructed from
the plumbing diagram specified by $G$, and the induced inner product
corresponds to the intersection form.  Note that the induced
intersection form on $V$ is negative-definite if and only if the
orbifold degree is negative.

A characteristic vector $K$ for $H_2(W(G))$ is a vector $K$ in the
dual space for $V$ with the property that $\langle K,v\rangle \equiv
m(v)\pmod{2}$ for each vertex $v$.  A sequence $\{K_i\}_{i=1}^\ell$ of
characteristic vectors is called a {\em full path} if
\begin{itemize}
\item for each $i$ and each vertex $v$ for $G$, 
$$|\langle K_i, v\rangle|\leq -m(v)$$
\item for each vertex $v$ for $G$, 
$$m(v)+2\leq \langle K_1,v \rangle \leq -m(v)$$
and
$$m(v)\leq \langle K_\ell,v \rangle \leq -m(v)-2$$
\item for each $i<\ell$, there is a vertex $v$ with the property that
$\langle K_i,v\rangle = -m(v)$, and $K_{i+1}=K_i+2\PD[v]$.
\end{itemize}
We call two full paths equivalent if they start with the same initial
vector.  As proved in~\cite{SomePlumbs}, equivalent full paths also
have the same final vector.

Full paths are related to the Heegaard Floer homology of Seifert fibered
spaces, according to the following result from~\cite{SomePlumbs}:

\begin{theorem} 
\label{thm:Seiferts}
Let $Y$ be a Seifert fibered rational
homology sphere with Seifert invariants 
$(b; \beta_1/\alpha_1, ... \beta_n/\alpha_n)$,
and let $G$ denote its corresponding negative-definite
plumbing graph. If $b\leq -n$, then $Y$ is an $L$-space.
More generally, $Y$ is an $L$-space if and only 
if the number of equivalence classes of full paths for $G$
agrees with the number of elements $|H_1(Y;\Z)|$.
\end{theorem}

\begin{proof}
The first statement is a consequence of  Proposition~\ref{prop:LSpaces}
and indeed, a proof is spelled out in~\cite{HolDiskSymp}.  
The second
is an application of the main result in~\cite{SomePlumbs}.
\end{proof}

\begin{prop}
\label{prop:SpaceForms}
Every three-manifold with elliptic geometry is an $L$-space.
\end{prop}

\begin{proof}
Spaces with elliptic geometry are those Seifert fibered fibered spaces
over a base orbifold $\Sigma$ with positive orbifold Euler
characteristic $\chi^{\orb}(\Sigma)$, which have non-zero orbifold
degree over their base (see~\cite{Scott} for a discussion of these
notions). 

Now, positivity of the Euler characteristic of the base forces it to
have genus zero and at most three singular fibers. If the number of
singular fibers is less than three, the total space is a lens space,
and hence covered by our earlier discussion. If there are three
singular fibers, with integral multiplicities $\alpha_1$, $\alpha_2$, and $\alpha_3$ (all
$>1$), then the formula for the orbifold Euler characteristic is
$$\chi^{\orb}(\Sigma)=-1+\frac{1}{\alpha_1}+\frac{1}{\alpha_2}+\frac{1}{\alpha_3}.$$ The
positivity criterion here forces
$\{\alpha_1,\alpha_2,\alpha_3\}$ to be either $\{2,3,n\}$ with $n\leq 5$
or $\{2,2,n\}$ with $n$ arbitrary.

By reversing the orientation of $Y$, we can arrange for the orbifold
degree to be negative.  
For these values of the $\alpha_i$,
negativity of
the  orbifold degree now
clearly forces $b\leq -2$. Indeed, by Theorem~\ref{thm:Seiferts}, it
suffices to consider the case where $b=-2$. 

We consider first the case where $\alpha_2=2$. In this case, we apply first
an induction on the length of the third chain of spheres, and a
subinduction on the multiplicity on the last leaf. For the basic case,
(where the graph has four vertices with multiplicities $-2$), we again
appeal to the second criterion from Theorem~\ref{thm:Seiferts}.  For
the inductive step, we consider the case where we add one more vertex
with multiplicity $-2$. It is easy to see that this is realized as
$Y_{\lambda+\mu}(K)$, where $Y$ is the Seifert fibered space obtained
by the graph where this last vertex is deleted, and $Y_{\lambda}(K)$
is obtained as a graph with fewer vertices (gotten by exchanging the
multiplicity with $-1$, and then successively blowing down $-1$
spheres).  Indeed, our hypotheses force
$$|H_1(Y_{\lambda+\mu}(K))|=|H_1(Y)|+|H_1(Y_{\lambda}(K))|,$$ so the
inductive hypotheses and Proposition~\ref{prop:LSpaces} applies to
show that $Y_{\lambda+\mu}(K)$ is an $L$-space. The induction required
to reduce the multiplicity of this last vertex by one works in the
same way.

Indeed, the case where $\alpha_2=3$ and $\beta_2=1$ follows from the
case where $\alpha_2=2$ by another application of
Proposition~\ref{prop:LSpaces}.

The finitely many (seven) cases
where $\beta_2=2$ and $\alpha_2=3$ which are not covered above
all follow from calculations using the
second criterion from Theorem~\ref{thm:Seiferts}.
\end{proof}

There are examples of non-elliptic Seifert fibered $L$-spaces.  In
fact, recall (see~\cite{BleilerHodgson}, see
also~\cite{FSLensSurgeries} and~\cite{Mattman}) that if we consider
the pretzel knot $P(-2,3,n)$ , where $n$ is an odd integer, then
$S^3_{2n+4}(K)$ is $\pm Y$ where $Y$ is the Seifert fibered rational
homology three-sphere with invariants $(-2; 1/2, 1/4, (n-8)/(n-6))$.  It
is easy to see that when $n\geq 9$, $S^3_{2n+4}(K)$ is an $L$-space
(though when $n>9$, it is not elliptic). As we mentioned in the
introduction, in the case where $n=7$, similar considerations show
that $S^3_{18}(K)$ is a lens space (c.f.~\cite{FSLensSurgeries}
and~\cite{BleilerHodgson}).  In view of
Proposition~\ref{prop:LSpaces}, if $n$ is any an odd integer with
$n\geq 7$, and $p$ be any integer with $p\geq 2n+4$,  then 
the
three-manifold $S^3_p(P(-2,3,n))$ is an $L$-space.

\section{Proof of Theorem~\ref{thm:FloerHomology}.}
\label{sec:Proof}

Theorem~\ref{thm:FloerHomology} follows from the relationship between
the knot Floer homology associated to $K\subset S^3$ and the Heegaard
Floer homology of $\HFa(S^3_p(K))$ for all large enough $p$. This
relationship is established in Section~\ref{Knots:sec:Relationship}
of~\cite{Knots} (see especially
Theorem~\ref{Knots:thm:LargePosSurgeries}). We recall these
constructions briefly here.

Recall that a knot $K$ induces a filtration on $\CFa(S^3)$.  More
precisely, the (finitely many) generators for $\CFa(S^3)$ have a
filtration level taking values in $(0,\Z)$ (the relevance of the first
coordinate will become apparent in a moment), and the differential is
non-increasing in this filtration. For $m\in\Z$, we let
$C\{(0,m)\}\subset \CFa(S^3)$ denote the subgroup generated by
elements with filtration level $(0,m)$.  More generally, we can let
$C\{(\ell,m)\}$ denote the set of generators of $C\{(0,m-\ell)\}$, now
shifted by a group isomorphism $$U^{\ell}\colon
C\{(\ell,m)\}\longrightarrow C\{(0,m-\ell)\}.$$

In~\cite{Knots}, we equip $C=\bigoplus_{(i,j)\in\Z}C\{(i,j)\}$
with a differential $D$ which commutes with the maps $U^\ell$,
and which is compatible with the initial differential on
$$\CFa(S^3)=\bigoplus_{m\in\Z} C\{(0,m)\}.$$ Indeed, the differential $D$ respects the
$\Z\oplus \Z$ filtration which sends the summand $C\{(i,j)\}\subset C$
to $(i,j)\in\Z\oplus\Z$. This means that if $\xi$ is supported in this
summand $C\{(\ell,m)\}$, then $\partial \xi$ is supported in the group 
$$\bigoplus_{\{(i,j)\big| i\leq \ell,~\text{and}~ j\leq m\}}
C\{(i,j)\}\subset C.$$
The complex $C$ is graded by the convention that 
$$H_*(C\{i=0\})\cong H_*(\CFa(S^3))\cong \Z$$
is supported in dimension zero, 
the differential $D$ lowers degree by one,
and the map  $U^{\ell}$ lowers it by $2\ell$. The chain complex
referred to here as $C$ is the $\Z\oplus\Z$-filtered complex $\CFKinf(S^3,K)$
from~\cite{Knots}.

A $\Z\oplus \Z$-filtered complex $C$ induces many other chain complexes.
We introduce the following notational shorthand.
If $R$ is a region in
the $(i,j)$ plane, then let $C(R)$ denotes the naturally induced complex
on the set of generators of $C$ whose filtration level $(i,j)$ lies in
the region $R$. Of course, this does not make sense for any region $R$
but there are three cases of interest to us (here, we write
$(i_1,j_1)\leq (i_2,j_2)$ if $i_1\leq i_2$ and $i_2\leq j_2$):
\begin{itemize}
\item suppose $R$
has the property that  
$$(i_1,j_1)\in R~\text{and}~(i_2,j_2)\leq (i_1,j_1)\Rightarrow
(i_2,j_2)\in R,$$ then  $C(R)$ is naturally a subcomplex;
\item suppose $R$ has the property that 
$$(i_1,j_1)\in R~\text{and}~(i_2,j_2)\geq (i_1,j_1)\Rightarrow
(i_2,j_2)\in R,$$ then  $C(R)$ is naturally a quotient complex;
\item suppose $R$ has the property that  
$$(i_1,j_1)\leq (i_2,j_2) \leq (i_3,j_3)~\text{and}~
(i_1,j_1), (i_3,j_3)\in R \Rightarrow (i_2,j_2)\in R,$$
then $C(R)$ is naturally  the subcomplex of a
quotient complex of $C$.
\end{itemize}

If $n$ is any integer, there is a natural affine identification
$\SpinC(S^3_n(K))\cong \Zmod{n}$ made explicit in~\cite{Knots} (but
not crucial for our present applications).  If $[m]\in\Zmod{n}$, we
let $\HFa(S^3,[m])$ denote the summand of the Floer homology in the
$\SpinC$ structure corresponding to $[m]$.
 Theorem~\ref{Knots:thm:LargePosSurgeries} of~\cite{Knots} states
that given any knot $K\subset S^3$,
there is an integer $N$ so that for all $n\geq N$,
there is an isomorphism of chain complexes
\begin{equation}
\label{eq:Identification}
\CFa(S^3_n(K),[m])\cong C\{\max(i,j-m)=0\}.
\end{equation}

Theorem~\ref{thm:FloerHomology} is now an algebraic consequence of
the above theorems. This algebra is encoded in the following two
lemmas.

\begin{lemma}
\label{lemma:EvenCase}
Let $C$ be a $\Z\oplus\Z$-filtered chain
complex over a field $\Field$, 
and let
$m$ be an integer with the property that
\begin{equation}
\label{eq:Hypotheses}
H_*(C\{\max(i,j-m)=0\})\cong \Field \cong H_*(C\{\max(i,j-(m-1))=0\})
\end{equation}
(ignoring the grading).

Suppose also that $H_*(C\{i<0,j=m\})=0$.

Then, either $H_*(C\{(0,m)\})=0$, in which case
$H_*(C\{i<0,j=m-1\})=0$ as well; or $H_*(C\{(0,m)\})\cong \Field$, in which
case $H_*(C\{i<0,j=m-1\})\cong \Field$ and  $H_*(C\{i=0,j\leq m-1\})=0$.
\end{lemma}

\begin{proof}
Let 
\begin{eqnarray*}
X=\{i\leq 0, j=m\}&{\text{and}}& Y=\{i=0,j\leq m-1\},
\end{eqnarray*}
so that $UX=\{i\leq -1,j=m-1\}$. In this notation,
Equation~\eqref{eq:Hypotheses} says that:
$$H_*(C\{UX\cup Y\})\cong H_*(C\{X\cup Y\})\cong \Field.$$

By the long exact sequence associated to the short exact sequence
$$
\begin{CD}
0@>>> C\{i<0,j=m\} @>>> C\{X\} @>>> C\{i=0,j=m\}@>>> 0,
\end{CD}
$$
combined with our hypothesis that $H_*(C\{i<0,j=m\})$ is trivial,
we see that 
$$H_*(C\{X\})\cong H_*(C\{(0,m)\}).$$

We have the following pair of short exact sequences
\begin{equation}
\label{eq:BigDiagram}
\begin{CD}
&	&	&	&	&		&	0	\\
&	&	&	&	&		&	@VVV	\\
0@>>> C\{U X\} @>{}>> C\{UX\cup Y \} @>{B}>> C\{Y\} @>>> 0	
\\
&	&	&	&	&		&	@V{}VV	\\
&	&	&	&	&		&	C\{X\cup Y\}	\\
&	&	&	&	&		&	@V{}VV	\\
&	&	&	&	&		&	C\{X\}	\\
&	&	&	&	&		&	@VVV	\\
&	&	&	&	&		&	0.
\end{CD}
\end{equation}

Let $n$ denote the rank of $H_*(C\{(0,m)\})\cong H_*(C\{X\})$.
There are two cases according to whether
or not the map induced by $B$ on homology
$$b\colon H_*(C\{UX\cup Y\}) \cong \Field \longrightarrow H_*(C\{Y\})$$
is trivial.

If $b$ is trivial, then it is easy to see that $H_*(C\{Y\})$ has
rank $n-1$, and that the coboundary associated to the horizontal short 
exact sequence
$$\delta_h\colon
H_*(C\{Y\})
\cong \Field^{n-1}\longrightarrow H_*(C\{UX\})\cong
\Field^{n}$$ 
is injective. Moreover, another count of ranks then ensures that the
coboundary map associated to the vertical short exact sequence
$$\delta_v\colon H_*(C\{X\}) \longrightarrow H_*(C\{Y\})$$ is
surjective. In particular, the image of the composite $$\delta_h\circ
\delta_v \colon H_*(C\{X\})\longrightarrow H_*(C\{UX\})$$ has rank
$n-1$. On the other hand, the relation that $D^2=0$ (on the induced
complex $C\{X\cup UX\cup Y\}$) ensures that the composite is trivial
(indeed, the relation gives a null-homotopy of the composite on the
chain level). Thus, we have established that $n=1$ and $H_*(C\{Y\})$ is 
trivial.

On the other hand, if $b$ is non-trivial, then
$H_*(C\{Y\})$ has rank $n+1$, and indeed the map
$$\delta_h\colon H_*(C\{Y\})\cong \Field^{n+1} \longrightarrow H_*(C\{UX\})\cong \Field^n$$
is surjective, while the map 
$$\delta_v\colon H_*(C\{X\})\cong \Field^{n} \longrightarrow H_*(C\{Y\})\cong \Field^{n+1}$$
is injective. 
On the one hand, this implies that the image of $(U\delta_v)\circ \delta_h$ is
$n$-dimensional; on the other hand, the composite is trivial (which
follows from the fact that $D^2=0$ on the complex $C\{UX\cup Y \cup UY\}$),
and hence $n=0$, and $H_*(C\{Y\})$ is one-dimensional.
These two cases cover the two cases in the conclusion of the lemma.
\end{proof}

\begin{lemma}
\label{lemma:OddCase}
Suppose once again that $C$ is a bigraded complex, with the property that
$$H_*(C\{\max(i,j-m+1)=0\})\cong \Field$$ for some integer $m$.

Suppose furthermore that $H_*(C\{i<0,j=m\})\cong \Field$ and
$H_*(C\{i=0,j\leq m\})=0$.
 
Then either $H_*(C\{(0,m)\})=0$, in which case
$H_*(C\{i<0,j=m-1\})=\Field$ and $H_*(C\{i=0,j\leq m-1\})=0$
as well; or $H_*(C\{i=0,j=m\})\cong \Field$, in which
case  $H_*(C\{i<0,m-1\})=0$.
\end{lemma}

\begin{proof}
We continue using the notation for the proof of Lemma~\ref{lemma:EvenCase},
with 
\begin{eqnarray*}
X=\{i\leq 0, j=m\}&{\text{and}}&
Y=\{i=0,j\leq m-1\}.
\end{eqnarray*} In this case, the hypothesis that
$$0=H_*(C\{i=0,j\leq m\})=H_*(C\{Y\cup(0,m)\})$$ ensures that
$$
H_*(C\{(0,m)\})\cong H_*(C\{Y\}).
$$
Let $n$ denote the rank of $H_*(C\{(0,m)\})$.

Again, we have the two exact sequences illustrated in the
Diagram~\eqref{eq:BigDiagram}.  Since $H_*(UX\cup Y)\cong \Field$ (and our
hypotheses also ensure the $H_*(X\cup Y)\cong \Field$), we have two cases
according to whether or not the map on homology $b$ trivial as before.

If $b$ is trivial, a diagram chase shows that
$$\delta_h\colon H_*(C\{Y\})\cong \Field^n\longrightarrow H_*(C\{UX\})\cong
\Field^{n+1}$$ is injective and also that $$\delta_v\colon
H_*(C\{X\})\longrightarrow H_*(C\{Y\})$$ is surjective. 
Again, this implies that the image of $\delta_h\circ \delta_v$
is $n$-dimensional, but since the composite is trivial, 
$n=0$,
and hence $H_*(C\{UX\})$ is one-dimensional.

If $b$ is non-trivial, 
$$\delta_h\colon H_*(C\{Y\})\cong \Field^n\longrightarrow H_*(C\{UX\})\cong
\Field^{n-1}$$ is surjective and also
$$\delta_v\colon H_*(C\{X\})\longrightarrow H_*(C\{Y\})$$ is injective.  
Thus, the image of $(U\delta_v)\circ \delta_h$ is
$(n-1)$-dimensional; but again this is trivial, forcing
$n=1$, and $H_*(C\{UX\})$ to be zero-dimensional. 

These two cases cover the two cases in the conclusion of the lemma.
\end{proof}

\vskip.2cm
\noindent{\bf{Proof of Theorem~\ref{thm:FloerHomology}.}}
By the universal coefficients theorem, it suffices to 
establish Theorem~\ref{thm:FloerHomology} over an arbitrary field
$\Field$. First, note that according to Equation~\eqref{eq:Identification}
(and our hypothesis on $K$), $H_*(C\{\{\max(i,j-m)\})=\Field$
for all $m$

Since $C\{i=0\}$ is finitely generated, for all sufficiently large
$m$, $C\{i<0,j=m\}=0$. Moreover, 
$C\{i=0,j\leq m\}=C\{i=0\}$, so that $H_*(C\{i=0,j\leq m\})\cong \Field$
is supported in even (zero) degree.
In particular,  the hypotheses of
Lemma~\ref{lemma:EvenCase} apply. 
Indeed, by descending induction on $m$, and using 
Lemma~\ref{lemma:EvenCase} and ~\ref{lemma:OddCase},
we have that for all $m$, the rank of $\HFKa(K,m)$ is 
at most one, and for each integer $m$, exactly one of
the following two possibilities  holds:
\begin{list}
{(\arabic{bean})}{\usecounter{bean}\setlength{\rightmargin}{\leftmargin}}
\item  either $H_*(C\{i<0,j=m\})=0$, 
	in the case where there either is no $\ell>m$ with $\HFKa(K,\ell)\neq 0$
	or the smallest such $\ell$ has the corresponding $\HFKa(K,\ell)$ supported
	in odd degree;
\item or $H_*(C\{i<0,j=m\})\cong \Field$ and $H_*(C\{i=0,j\leq m\})=0$ and the smallest
	$\ell>m$ with $\HFKa(K,\ell)\neq 0$ has the corresponding $\HFKa(K,\ell)$
	supported in even degree.
\end{list}

Indeed, let $\ell>m$ be a pair of integers for which
$\HFKa(K,\ell)\cong \Field\cong \HFKa(K,m)$ (ignoring gradings), and for 
all intermediate  $m<j<\ell$, $\HFKa(K,j)=0$. Let $d$ denote
the dimension in which $\HFKa(K,\ell)$ is supported. 

If $d$ is even, then it is easy to see
(once again, by one application
of Lemma~\ref{lemma:EvenCase} followed by repeated
applications of Lemma~\ref{lemma:OddCase}) that
$H_*(C\{i<0,j=m\})\cong \Field$ is supported in dimension $d-2(\ell-m)$.
It follows now from Lemma~\ref{lemma:OddCase}, that the coboundary
map for the short exact sequence
$$
\begin{CD}
0@>>>C\{i<0,j=m\} @>>> C\{i\leq 0,j=m\} @>>> \CFKa(K,m)@>>> 0,
\end{CD}
$$
which drops dimension by one, induces an isomorphism in homology
$$\delta_h\colon \HFKa(K,m)\longrightarrow H_*(C\{i<0,j=m\});$$
thus  $\HFKa(K,m)$ is supported in dimension $d-2(\ell-m)+1$.

If $d$ is odd, then it follows that the coboundary map for
the short exact sequence
$$
\begin{CD}
0@>>> C\{i=0,j<\ell \} @>>> C\{i=0,j\leq \ell\} @>>> \CFKa(K,\ell)@>>>0$$
\end{CD}
$$ induces an isomorphism on homology. Thus, $H_*(C\{i=0,j<\ell
\})\cong \Field$ is supported in dimension $d-1$. Indeed, it is
easy to see that the natural inclusion $C\{i=0,j\leq m\}\subset
C\{i=0,j<\ell \}$ induces an isomorphism in homology.  In fact since
$H_*(C\{i=0,j\leq m-1)=0$ (c.f. Lemma~\ref{lemma:EvenCase}), the
projection $C(\{i=0,j\leq m\})\longrightarrow C\{i=0,j=m\}$ also
induces a (degree-preserving) isomorphism in homology: i.e.
$\HFKa(K,m)$ is supported in dimension $d-1$.

Together, these claims establish the theorem.
\qed
\vskip.2cm

It is worth pointing out that the above lemmas hold even in the case
where $S^3_p(K)$ is not an $L$-space. In particular, we have the
following:

\begin{prop}
\label{prop:TauAlone}
Let $K\subset S^3$ be a knot with the properties that $\HFKa(K,m)=0$
for all $m>d$, $\HFKa(K,d)\neq 0$. Then, if for all sufficiently large
$n$, $\HFa(S^3_n(K),[d])\cong \HFa(S^3_n(K),[d-1])\cong \Q$, we have that
$\tau(K)=d$.
\end{prop}

\begin{proof}
Follows immediately from Lemma~\ref{lemma:EvenCase}.
\end{proof}

We turn now to some of the consequences of
Theorem~\ref{thm:FloerHomology} which were described in the introduction.

\vskip.2cm
\noindent{\bf{Proof of Corollary~\ref{cor:StructAlex}.}}
This is an immediate consequence of Theorem~\ref{thm:FloerHomology}
and the relationship between the knot Floer homology and
the Alexander polynomial, 
Proposition~\ref{Knots:prop:Euler} of~\cite{Knots}.
\qed

\vskip.2cm
\noindent{\bf{Proof of Corollary~\ref{cor:LensCondition}.}}
Note that the above proof could be modified to give the relationship
between the absolute gradings on $\HFa(L)$ with the Alexander
polynomial of $K$. We have not spelled this out, as it 
was already determined in~\cite{AbsGraded}, c.f. 
Theorem~\ref{AbsGraded:thm:PSurgeryLens} and especially 
Corollary~\ref{AbsGraded:cor:AlexLens}, both in~\cite{AbsGraded}.
\qed

\vskip.2cm
\noindent{\bf{Proof of Proposition~\ref{prop:Experiment}.}}
Proposition~\ref{prop:Experiment} is verified by first calculating
$d(-L(p,q),i)$ for $p$ in some range, then enumerating all possible
correspondences $\sigma$, keeping only those $q$ for which one of the
correspondences $\sigma$ satisfies the conditions of
Corollary~\ref{cor:LensCondition}, and then verifying that this list
of allowed lens spaces is covered by Berge's list.  This verification
is algorithmic, if tedious. We used code written in
Mathematica~\cite{Mathematica}.
\qed
\section{Alternating knots and $L$-space surgeries}
\label{sec:AltKnots}

For the proof of Theorem~\ref{thm:AltLens}, we use the following
characterization of the $(2,n)$ torus knots, which follows easily from
standard properties of the Alexander polynomial for alternating knots,
compare~\cite{MurasugiAlt}, \cite{Crowell}, ~\cite{Kauffman},
~\cite{Menasco}:

\begin{prop}
\label{prop:AlexAlt}
If $K$ is an alternating knot with the property that all
the coefficients $a_i$ of its Alexander polynomial $\Delta_K$ 
have $|a_i|\leq 1$, then $K$ is the $(2,2n+1)$ torus knot.
\end{prop}

\begin{proof}
According to a theorem of Menasco~\cite{Menasco}, a non-prime
alternating  knot factors as a sum of (non-trivial)
alternating knots. According to a theorem of Crowell and
Murasugi~\cite{Crowell} and~\cite{MurasugiAlt}, the Alexander
polynomials of these factors are non-trivial polynomials whose
coefficients alternate in sign. It follows at once that the Alexander
polynomial of our original knot has coefficients greater than one.

Thus, it suffices to consider the case where $K$ is prime.  Consider
an alternating projection, and let $w$ resp. $b$ denote the number of
white resp. black regions in the checkerboard coloring. Form the
``black graph'' $B$ of the knot projection, whose vertices correspond
to the black regions and whose edges correspond to double-points in
the knot projection.  Recall that the Alexander polynomial of a knot
can be interpreted as a suitable count of spanning trees of $B$, where
each tree is weighted by some $T$-power, and a sign.  Moreover, the
number of distinct $T$-powers appearing this polynomial is bounded
above by the number double-points in the knot projection plus one, which in
turn is given by $w+b-1$.  Recall that the main step in the
Crowell-Murasugi theorem shows that for an alternating projection, the
trees contributing a fixed $T$-power contribute with the same sign.

According to a result of Crowell~\cite{CrowellTwo}, the total number
of such trees for a prime, alternating knot is bounded below by
$1+(w-1)(b-1)$. Thus, according to our hypothesis that $|a_i|\leq 1$,
no two trees can contribute to the same $T$-power, and hence
$$1+(w-1)(b-1)\leq w+b-1.$$ This inequality immediately forces either
$w=b=3$ or at least one of $w$ or $b=2$. In the case where $w=b=3$, it
is easy to see that the knot in question is the figure eight knot,
whose Alexander polynomial does not satisfy the hypotheses of the
theorem. In the case where $w$ or $b=2$, it is easy to see that $K$ is
the $(2,2n+1)$ torus knot.
\end{proof}

\vskip.2cm
\noindent{\bf{Proof of Theorem~\ref{thm:AltLens}.}}
Put together Corollary~\ref{cor:StructAlex} and 
Proposition~\ref{prop:AlexAlt}.
\section{Berge knots are fibered}
\label{sec:BergeKnots}

The aim of this section is to verify
Proposition~\ref{prop:BergesFibered}.  This result seems to be known
to the experts, but we include a proof here for completeness.

According to a theorem of Stallings~\cite{Stallings}, a connected
three-manifold $Y$ with $H_1(Y;\Z)\cong \Z$ is fibered if and only if
the kernel of the Abelianization map $$\pi_1(Y)\longrightarrow
H_1(Y;\Z)$$ is a finitely generated group.  We use this
characterization to show that all knots coming from Berge's
construction (c.f. Definition~\ref{def:BergeKnots}) are fibered.
Specifically, the knot complements arising from Berge's constructions
have Heegaard genus two, and thus their fundamental group admits a
presentation with two generators and one relator $$G=\langle X,
Y\rangle/R(X,Y).$$ A theorem of Brown~\cite{Brown} concerns conditions
under which a homomorphism $$\chi\colon G
\longrightarrow \R$$ has finitely generated kernel. Specifically, 
write the word $$R(X,Y)= A_1\cdot ... \cdot A_m,$$ where $A_i\in
\{X,Y,X^{-1},Y^{-1}\}$, and then consider the sequence of real numbers
$$S=\left\{\chi\left(\prod_{i=1}^n A_i\right) \right\}_{n=1}^m,$$ then Brown's
theorem states that the kernel of $\chi$ is finitely generated if the
sequence $S$ achieves its maximum and minimum only once. We will apply
this condition to the Abelianization map in the following proof of
Proposition~\ref{prop:BergesFibered}:

\vskip.3cm
\noindent{\bf{Proof of Proposition~\ref{prop:BergesFibered}.}}
Berge's construction gives a knot for each generator for the homology
of $H_1(L(p,q);\Z)$, also giving rise to a genus two Heegaard diagram
for the knot complement $L(p,q)-K$. Explicitly, we start with a genus
one Heegaard diagram for $L(p,q)$: the Heegaard surface is given as a
square torus, and $\alpha$ is a straight line with slope $p/q$, while
$\beta_1$ is given as a line with slope zero. Now, fix an integer $0<
k <p$ which is relatively prime to $p$, and draw a segment with slope
zero which is disjoint from $\beta_1$ and which intersects $\alpha$ in
$k$ points. After attaching a one-handle to the torus at near the
endpoints of the arc, we can close up the arc to give a closed circle
$\beta_2$ in the surface of genus two which continues to meet $\alpha$
in $k$ points and is disjoint from $\beta_1$. The associated Heegaard
diagram is easily seen to represent a knot complement $L(p,q)-K$,
where $K$ is a lens space Berge knot in the sense of
Definition~\ref{def:BergeKnots} representing a homology class which is
$k$ times a generator for $H_1(L(p,q),\Z)$.

\begin{figure}
\mbox{\vbox{\epsfbox{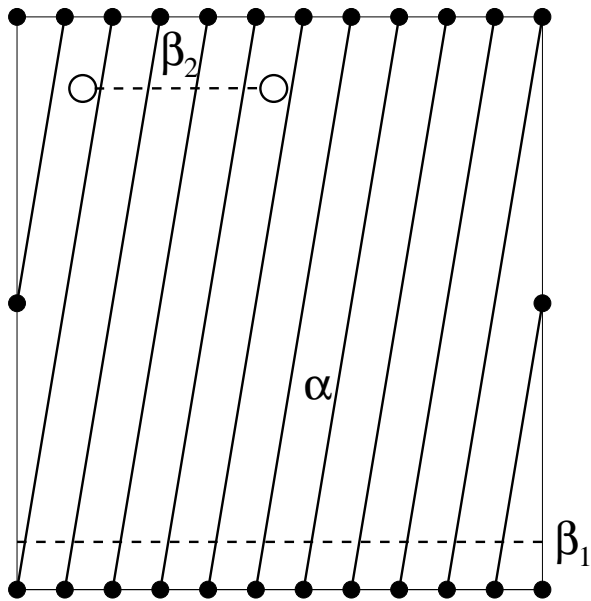}}}
\caption{\label{fig:ThreeFour}
{\bf{The $(3,4)$ torus knot.}} We draw here the Heegaard diagram
for a Berge knot in $L(11,2)$. This drawing takes place in a square torus
(i.e. make the usual identifications on this square), with an
additional one-handle added along the two hollow circles. The
$\alpha$-curve is the diagonal curve with slope $11/2$, $\beta_1$ is
the long horizontal dashed line, and $\beta_2$ has an arc indicated by
the other dashed line, which then closes up inside the attached handle.
Tracing along $\alpha$, we see that  the fundamental group of
$L(11,2)-K$ is generated by elements
$X$ and $Y$
satisfying the relation $XYXYX^5YXYX^3=e$.
Indeed, the Heegaard diagram we obtain in this manner describes
the complement of the $(3,4)$ torus knot $T_{3,4}$ in $S^3$
(corresponding to the fact that $+11$ surgery on $T_{3,4}$ gives
$-L(11,2)$).}
\end{figure}

This description can be used to give a presentation of the fundamental
group of $Y=L(p,q)-K$. Of course, the description gives $Y$ as a genus
two handlebody and an attached disk, and hence $\pi_1(Y)$ as a group
with two generators and one relation. This description can be given
explicitly: let $X$ and $Y$ be the curves dual to the attaching disks
for $\beta_1$ and $\beta_2$. The relation arising from the attaching
disk $\alpha$ is found by following the curve $\alpha$, and recording
in order which of the curves $\beta_1$ and $\beta_2$ are encountered
-- with $\beta_1$ contributing a factor of $X^\pm$ and $\beta_2$
contributing a factor of $Y^\pm$, where here the exponent is given by
the local intersection number of $\alpha$ with the corresponding
$\beta$-curve. In fact, for the curves coming from Berge's construction,
the obtained relator has the following simple form.
Let $$E(i)=E(i,p,q,k)=
\left\{
\begin{array}{ll}
1 &{\text{if there is some integer $0\leq j <k$ with
		$j\equiv i\cm q \pmod{p}$}} \\
0 &{\text{otherwise,}}
\end{array}
\right.$$
then the presentation of $G=\pi_1(L(p,q)-K)$
given by the above procedure is
$$G\cong  \langle X, Y\rangle /\Pi_{i=1}^p
\left(X Y^{E(i,p,q,k)}\right).$$  
See Figure~\ref{fig:ThreeFour} for an example.

It follows that $G/[G,G]$ is the lattice spanned by $[X]$ and $[Y]$,
modulo the relation $p[X]+k[Y]=0$. Thus, Abelianization can be viewed
as a map $$\chi\colon G \longrightarrow \Z$$ which sends $[X]$ to $-k$
and $[Y]$ to $p$. 

Now, our relator
$R(X,Y)=\prod_{i=1}^p(XY^{E(i,p,q,k)})$ contains $k$ instances of $Y$;
explicitly, writing $R(X,Y)=\prod_{i=1}^m A_i$, there is a sequence of
distinct integers $\{n_i\}_{i=1}^k$ with the property that
$A_{n_i}=Y$. It is easy to see that the maxima of the
sequence $S$ described above are achieved amongst the $k$ words of the
form $\{w_i=A_1\cm...\cm A_{n_i}\}_{i=1}^k$.  Moreover, it is
straightforward to see that 
$$\chi(w_i)\equiv i\cm p \pmod{k},$$ and
hence, since $(k,p)=1$, these $k$ values are distinct, showing
uniqueness of the maximum.
Similarly, the minima are achieved on the $k$ words
$\{u_i=A_1\cm...\cm A_{n_i-1}\}$.
Once again, these $k$ values $\chi(u_i)=\chi(w_i)-p$
are distinct modulo $k$, and hence
the minimum is uniquely achieved. It follows now from Brown's
theorem~\cite{Brown} that the kernel of the Abelianization map is
finitely generated, and hence according to Stallings'
theorem~\cite{Stallings} that the knot complement is fibered.
\qed

\commentable{
\bibliographystyle{plain}
\bibliography{biblio}
}

\end{document}